\documentclass[graybox]{svmult}

\usepackage{mathptmx}     
\usepackage{helvet}      
\usepackage{courier}   
\usepackage{type1cm}  
\usepackage{multicol}    
\usepackage[bottom]{footmisc}
\usepackage{makeidx,epic,eepic,amsfonts,color,latexsym,mathrsfs,amssymb,amsmath,epsfig,rotating,graphicx}

\makeindex 

\newtheorem{thm}{Theorem}[section]

\newtheorem{dfn}[thm]{Definition}
\newtheorem{algor}[thm]{Algorithm}
\newtheorem{cor}[thm]{Corollary} 
\newtheorem{rem}[thm]{Remark}
\newtheorem{prop}[thm]{Proposition}
\newtheorem{lem}[thm]{Lemma}
\newtheorem{ex}[thm]{Example}
\newlength{\smale}
\newlength{\jmr}
\newlength{\khov}
\newlength{\bernd}
\makeatletter
\renewcommand*\env@matrix[1][c]{\hskip -\arraycolsep
  \let\@ifnextchar\new@ifnextchar
  \array{*\c@MaxMatrixCols #1}}
\makeatother

\newcommand{\thth}{^{\text{\underline{th}}}}

\newcommand{\stst}{^{\text{\underline{st}}}}

\newcommand{\np}{{\mathbf{NP}}}

\newcommand{\Log}{\mathrm{Log}}

\newcommand{\Q}{\mathbb{Q}}
\newcommand{\R}{\mathbb{R}}

\newcommand{\C}{\mathbb{C}}
\newcommand{\N}{\mathbb{N}}
\newcommand{\Z}{\mathbb{Z}}
\newcommand{\bO}{\mathbf{O}}

\newcommand{\cH}{\mathcal{H}}
\newcommand{\cL}{\mathcal{L}}
\newcommand{\1}{{\text{\rm 1} \kern -.41em \text{\rm 1} }}

\newcommand{\Zn}{\Z^n}

\newcommand{\Qn}{\Q^n}
\newcommand{\Rn}{\R^n}

\newcommand{\Cn}{\C^n}

\newcommand{\Cs}{\C^*}

\newcommand{\Csn}{{(\C^*)}^n}

\renewcommand{\qed}{$\blacksquare$}
\newcommand{\dia}{$\diamond$}
\newcommand{\amoeba}{\mathrm{Amoeba}}
\newcommand{\newt}{\mathrm{Newt}}
\newcommand{\anewt}{\mathrm{ArchNewt}}
\newcommand{\trop}{\mathrm{ArchTrop}}

\newcommand{\conv}{\mathrm{Conv}}
\newcommand{\size}{\mathrm{size}}

\newcommand{\sign}{\mathrm{sign}}

\author{Eleanor Anthony \and Sheridan Grant \and Peter Gritzmann \and 
J.\ Maurice Rojas}
\institute{Eleanor Anthony \at Department of Mathematics, University of 
Mississippi, Hume Hall 305,
P.\ O.\ Box 1848, University, Mississippi \ 38677-1848, 
\email{ecanthon@go.olemiss.edu} . Partially supported by NSF REU grant 
DMS-1156589. \and 
Sheridan Grant \at Department of Mathematics, 
640 North College Avenue, 
Claremont, Calif.\ 91711, \email{sheridan.grant@pomona.edu} . 
Partially supported by NSF REU grant
DMS-1156589. \and  
Peter Gritzmann \at Fakult\"at f\"ur Mathematik, Technische Universit\"at 
M\"unchen, D-80290 M\"unchen, Germany, \email{gritzmann@tum.de} . Work 
supported in part by the German Research Foundation (DFG). \and 
J.\ Maurice Rojas \at Department of Mathematics,
Texas A\&M University TAMU 3368,
College Station, Texas \ 77843-3368, USA, \email{rojas@math.tamu.edu} . 
Partially supported by NSF MCS grant DMS-0915245 and Sandia National 
Laboratories. } 

\begin{document}

\title*{Polynomial-Time Amoeba Neighborhood Membership and Faster Localized 
Solving}  
\titlerunning{Faster Amoeba Neighborhoods and Faster Solving} 

\maketitle 

\date{\today} 

\abstract{  
We derive efficient algorithms for coarse approximation of algebraic 
hypersurfaces, useful for estimating the distance between an input 
polynomial zero set and a given query point. Our methods work best on sparse 
polynomials of high degree 
(in any number of variables) but are nevertheless completely general. 
The underlying ideas, which we take the time to describe in an elementary way, 
come from tropical geometry. We thus  
reduce a hard algebraic problem to high-precision linear optimization, 
proving new upper and lower complexity estimates along the way. }  

\mbox{}\\ 
\mbox{}\hfill {\em Dedicated to Tien-Yien Li, in honor of his birthday.} 
\hfill\mbox{} 

\section{Introduction} 
As students, we are often asked to draw, hopefully without 
a calculator, real zero sets of low degree polynomials in 
few variables. As scientists and engineers, we are often asked  
to count or approximate, hopefully with some computational assistance, 
real and complex solutions of arbitrary systems of polynomial equations in 
many variables. If one allows sufficiently coarse approximations, then the 
latter problem is as easy as the former. Our main results clarify 
this transition from hardness to easiness. In particular, we significantly 
speed up certain queries involving distances 
between points and algebraic hypersurfaces 
(see Theorems \ref{thm:dist}--\ref{thm:mixed} and 
Remark \ref{rem:speed} below). 

Polynomial equations are ubiquitous in numerous applications, 
such as algebraic statistics \cite{mle}, chemical reaction kinetics 
\cite{signdick}, discretization of partial differential equations 
\cite{hautumor}, satellite orbit design \cite{mortari}, 
circuit complexity \cite{kpr}, and 
cryptography \cite{hfe}. The need to solve larger and larger 
equations, in applications as well as for theoretical purposes, 
has helped shape algebraic geometry and numerical analysis for centuries. 
More recent work in algebraic complexity tells us that many basic questions 
involving polynomial equations are $\np$-hard (see, e.g., \cite{plaisted,
hnam,cxhilb,cxcountcompo}). This is by no means an excuse to consider 
polynomial equation solving hopeless: computational 
scientists solve problems of near-exponential complexity every day. 

More to the point, thanks to recent work on {\em Smale's 17th Problem} 
\cite{bps17,bcs17}, we have learned that randomization and approximation are 
the key to avoiding the bottlenecks present in hard deterministic questions 
involving roots of polynomial systems. Smale's 17th Problem concerns the 
complexity of approximating a {\em single} complex root of a random polynomial 
system and is well-discussed in \cite{21a,21b,ss1,ss2,ss3,ss4,ss5}. Our 
ultimate goal is to extend this philosophy to the harder problem 
of {\em localized solving}: estimating how far the {\em nearest} root of a 
given system of polynomials (or intersection of several zero sets) is from a 
given point. We make some initial steps by first approximating the shape of a 
single zero set, and we then outline a tropical-geometric approach to 
localized solving in Section \ref{sec:exp}. 

Toward this end, let us first recall the natural idea \cite{virologpaper} of 
drawing zero sets on log-paper. In what follows, 
we let $\C^*$ denote the non-zero complex numbers and write  
$\C\!\left[x^{\pm 1}_1,\ldots,x^{\pm 1}_n\right]$ for the ring of 
Laurent polynomials with complex coefficients, i.e., polynomials with 
negative exponents allowed. 

\noindent
\begin{picture}(200,200)(3,-120)
\put(312,-103){
\begin{minipage}[b]{0.3\linewidth}
\vspace{0pt}
\epsfig{file=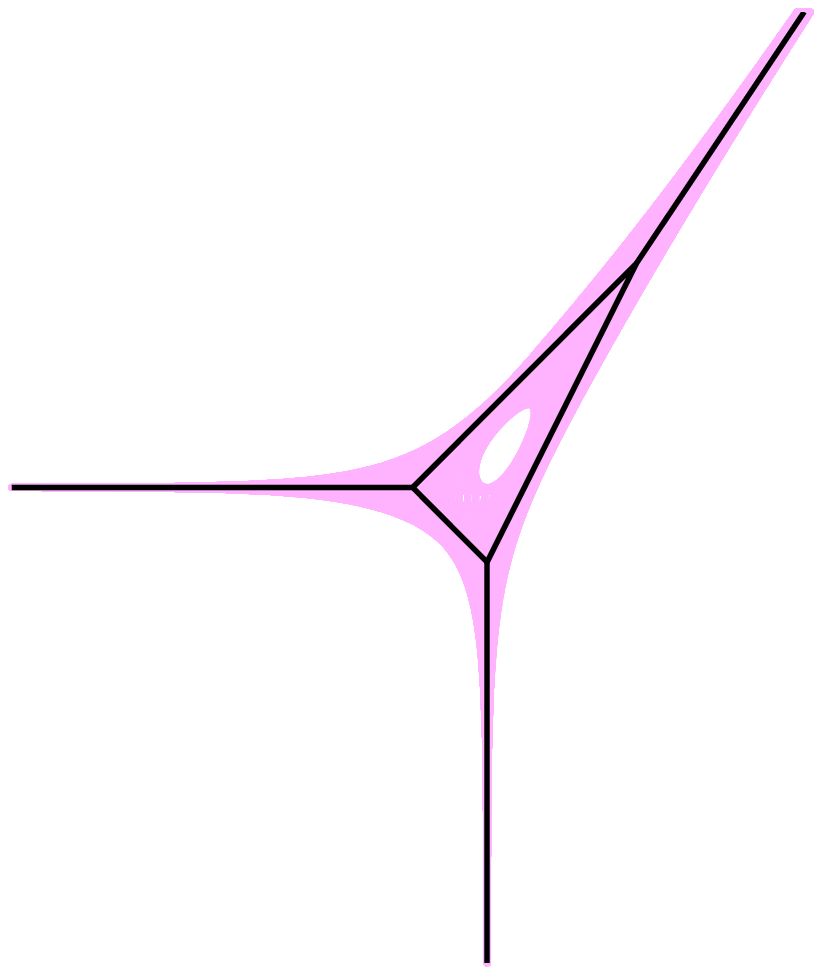,height=2.7in,clip=}
\end{minipage}}
\put(0,-100){
\begin{minipage}[b]{1\linewidth}
\vspace{0pt}
\begin{dfn}
\label{dfn:amoebaarchtrop}
\scalebox{1}[1]{We use the abbreviations $x\!:=\!(x_1,\ldots,x_n)$ 
and $\Log|x|\!:=\!(\log|x_1|, \ldots,\log|x_n|)$, and, for any}\\ 
\scalebox{1}[1]{$f\!\in\!\C\!\left[x^{\pm 1}_1,\ldots,x^{\pm 1}_n\right]$, we
define $\amoeba(f)\!:=\!\{\Log|x| \; | \; f(x)\!=\!0 \ , \ x\!\in\!\Csn\}$. 
We call $f$ an {\em $n$-variate}}\\
\scalebox{.95}[1]{{\em $t$-nomial} when we 
can write $f(x)\!=\!\sum^t_{i=1} c_i x^{a_i}$ with $c_i\!\neq\!0$, 
$a_i\!:=(a_{1,i},\ldots,a_{n,i})$, 
and $x^{a_i}\!:=\!x^{a_{1,i}}_1 x^{a_{2,i}}_2 \cdots x^{a_{n,i}}_n$}\\   
\scalebox{1}[1]{for all $i$. Finally, we define the {\em Archimedean tropical 
variety of $f$}, denoted $\trop(f)$, to be}\\ 
\scalebox{1}[1]{the set of all $w\!\in\!\Rn$ for
which $\max_i \left|c_ie^{a_i\cdot w}\right|$ is attained for at least 
{\em two} distinct indices $i$.$^1$ \dia }\\ 
\end{dfn}
\begin{ex} \label{ex:archtrop}
\scalebox{.98}[1]{Taking $f(x)\!=\!1+x^3_1+x^2_2-3x_1x_2$,
an illustration of $\amoeba(f)$ and $\trop(f)$,}\\ 
\scalebox{1}[1]{truncated to $[-7,7]^2$, appears to the  
right. $\amoeba(f)$ is lightly shaded, while $\trop(f)$ is}\\  
the piecewise-linear curve. \dia 
\end{ex} 
\scalebox{.99}[1]{One may be surprised that $\amoeba(f)$ and $\trop(f)$ are
so highly structured: $\amoeba(f)$ has}\\ 
tentacles reminiscent of a living amoeba, and $\trop(f)$ is 
a polyhedral complex, i.e., a union\\ 
of polyhedra intersecting only along common
faces. One may also be surprised that $\amoeba(f)$\\ 
\scalebox{.99}[1]{and $\trop(f)$ are so
closely related: for our example above, one set is strictly contained in the}  
\end{minipage}}
\end{picture} \addtocounter{footnote}{1} 
\footnotetext{Throughout this 
paper, for any two vectors $u\!:=\!(u_1,\ldots,u_N)$ and
$v\!:=\!(v_1,\ldots,v_N)$ in $\R^N$, we use $u\cdot v$ to denote the standard 
dot product $u_1v_1+\cdots+u_Nv_N$.} 

\vspace{-.7cm} 
\noindent 
other, every point of one set is close to some point of the other, and both
sets have topologically\\ 
similar complements (4 open connected components,
exactly one of which is bounded).

Proving that $\amoeba(f)$ and $\trop(f)$ are in fact {\em equal} 
when $f$ has two or fewer monomial terms is a simple exercise 
(see Proposition \ref{prop:binom} below).  
More generally, to quantify exactly how close $\amoeba(f)$ and $\trop(f)$ are,  
one can recall the {\em Hausdorff distance}, denoted $\Delta(U,V)$, between 
two subsets $U,V\!\subseteq\!\Rn$: it is defined to be the maximum 
of $\sup_{u\in U}{} \inf_{\substack{\mbox{}\\ v\in V}}|u-v|$ and 
$\sup_{v\in V}{} \inf_{\substack{\mbox{}\\ u\in U}}|u-v|$. 
We then have the following recent result of Avenda\~no, Kogan, Nisse, and 
Rojas. 
\begin{thm}
\label{thm:metric} 
\cite{aknr} 
For any $n$-variate $t$-nomial $f$ 
we have  
$\displaystyle{\Delta(\amoeba(f),\trop(f))\leq (2t-3)
\log(t-1)}$. 
In particular, we also have 
$\displaystyle{\sup
\limits_{\text{\scalebox{.7}[1]{$u\in\amoeba(f)$}}} \inf
\limits_{\substack{\mbox{} \\ 
\text{\scalebox{.7}[1]{$v\in\trop(f)$}}}} |u-v|\leq\log(t-1)}$. Finally, 
for any $t\!>\!n\!\geq\!1$, there is an $n$-variate 
$t$-nomial $f$ with $\Delta(\amoeba(f),\trop(f))\!\geq\!\log(t-1)$. \qed  
\end{thm}

\noindent 
\scalebox{.96}[1]{Note that the preceding upper bounds are completely 
independent of the coefficients, degree, and number of variables of $f$}
\linebreak 
We conjecture that an $O(\log t)$ 
upper bound on the above Hausdorff distance is possible. More practically, 
as we will see in later examples, $\amoeba(f)$ and $\trop(f)$ 
are often much closer than guaranteed by any proven upper bound. 

Given the current state of numerical algebraic geometry and algorithmic 
polyhedral geometry, the preceding metric result suggests that it might be 
useful to apply Archimedean tropical varieties to speed up polynomial system 
solving. Our first two main results help set the stage for such speed-ups. 
Recall that $\Q[\sqrt{-1}]$ denotes those complex numbers whose 
real and imaginary parts are both rational. 
Our complexity results will all be stated relative to the 
classical Turing (bit) model, with the underlying notion of input size 
clarified below in Definition \ref{dfn:input}.  
\begin{thm} 
\label{thm:dist} 
Suppose $f\!\in\!\C\!\left[x^{\pm 1}_1,\ldots,x^{\pm 1}_n\right]$ 
and $w\!\in\!\Rn$. Then\\ 
\mbox{}\hfill $-\log(t-1)\!\leq\!\Delta(\amoeba(f),w)-\Delta(\trop(f),w)
\!\leq\!(2t-3)\log(t-1)$.\hfill\mbox{}\\ 
In particular, if we also assume that 
$n$ is fixed and $(f,w)\!\in\!\Q[\sqrt{-1}]
\left[x^{\pm 1}_1,\ldots,x^{\pm 1}_n\right]
\times \Qn$ with $f$ a $t$-nomial, then we can compute polynomially many bits 
of $\Delta(\trop(f),w)$ in polynomial-time, and 
there is a polynomial-time\linebreak algorithm that 
declares either (a) $\Delta(\amoeba(f),w)\!\leq\!(2t-2)\log(t-1)$ 
or\\ 
\mbox{}\hspace{1.67in}(b) $w\!\not\in\!\amoeba(f)$ and 
$\Delta(\amoeba(f),w)\!\geq\!\Delta(\trop(f),w)-\log(t-1)\!>\!0$. 
\end{thm}
\begin{thm} 
\label{thm:region} 
Suppose $n$ is fixed. Then there is a polynomial-time algorithm that, for any 
input\linebreak 
$(f,w)\!\in\!\Q[\sqrt{-1}]\!\left[x^{\pm 1}_1,\ldots,
x^{\pm 1}_n\right]\times \Qn$ with $f$ a $t$-nomial, outputs the closure of 
the unique cell $\sigma_w$ of $\Rn\!\setminus\!\trop(f)$ (or $\trop(f)$) 
containing $w$, described as an explicit intersection of $O(t^2)$ half-spaces. 
\end{thm}  

\noindent 
The importance of Theorem \ref{thm:dist} 
is that deciding whether an input point $w$ lies in an input $\amoeba(f)$, 
{\em even restricting to the special case $n\!=\!1$}, is already $\np$-hard 
\cite{aknr}. Theorem \ref{thm:region} enables us to find explicit regions, 
containing a given query point $w$, where $f$ can not vanish. As we will 
see later in Sections \ref{sub:lo} and \ref{sub:log}, improving Theorems 
\ref{thm:dist} and \ref{thm:region} to 
{\em polynomial} dependence in $n$ leads us to deep questions in 
Diophantine approximation and the 
complexity of linear optimization. 

It is thus natural to speculate that tropical 
varieties can be useful for localized polynomial system solving, i.e., 
estimating how far the nearest root of a given system of $n$-variate 
polynomials $f_1,\ldots,f_k$ is from an input point $x\!\in\!\C$. Our 
framework indeed enables new positive and negative results on this problem. 
\begin{thm} 
\label{thm:mixed} 
Suppose $n$ is fixed. Then there is a polynomial-time algorithm that, for any 
input $k$ and\linebreak 
$(f_1,\ldots,f_k,w)\!\in\!\left(\Q[\sqrt{-1}]\!\left[x^{\pm 1}_1,
\ldots,x^{\pm 1}_n\right]\right)^k\times \Qn$,
outputs the closure of the unique cell 
$\sigma_w$ of $\Rn\setminus\bigcup^k_{i=1} \trop(f_i)$ (or 
$\displaystyle{\bigcap\limits_{\trop(f_i)\ni w} \!\!\!\!\!\!\!\!\!\! 
\trop(f_i)}$) containing $w$, described as an explicit intersection of 
half-spaces.  However, if $n$ is allowed to\linebreak 

\vspace{-.4cm}
\noindent
vary, then deciding whether $\sigma_w$ has a vertex in $\displaystyle{
\bigcap\limits^n_{i=1}\trop(f_i)}$ is 
$\np$-hard. 
\end{thm} 

\noindent 
We will see in Section \ref{sec:exp} how the first assertion is 
useful for finding special start-points for {\em Newton Iteration} 
and\linebreak 
{\em Homotopy Continuation} that sometimes enable the approximation of just 
the roots with norm vector near $(e^{w_1},\ldots,e^{w_n})$. The second 
assertion can be 
considered as a refined tropical analogue to a classical algebraic complexity 
result: deciding whether an arbitrary input system of polynomials equations 
(with integer coefficients) has a complex root is $\np$-hard \cite{gj}. 
However, in light of the recent partial solutions to Smale's 17th Problem  
\cite{bps17,bcs17} (showing that randomization and approximation help us  
evade $\np$-hardness for {\em average-case} inputs), we suspect that an 
analogous speed-up is possible in the tropical case as well. 

On the practical side, we point out that the algorithms underlying 
Theorems \ref{thm:dist}--\ref{thm:mixed} are quite easily 
implementable. (A preliminary {\tt Matlab} 
implementation of our algorithms is available upon request.) 
Initial experiments, discussed in Section \ref{sec:exp} below, indicate that a 
large-scale implementation could be a worthwhile companion to 
existing polynomial system solving software. 

Theorems \ref{thm:dist}, \ref{thm:region}, and \ref{thm:mixed} are 
respectively proved in Sections \ref{sec:dist}, \ref{sec:region}, and 
\ref{sec:mixed}. Before moving on to the necessary technical background, let 
us first clarify our underlying input size and point out some historical 
context. 
\begin{dfn} 
\label{dfn:input} 
We define the {\em input size} of a polynomial $f\!\in\!\Z[x_1,\ldots,x_n]$, 
written $f(x)\!=\!\sum^t_{i=1} c_i x^{a_i}$, to be 
$\size(f)\!:=\!\sum^t_{i=1} \log\left((2+|c_i|)\prod^n_{j=1}(2+|a_{i,j}|)
\right)$, where $a_i\!=\!(a_{i,1},\ldots,a_{i,n})$ for all $i$.
Similarly, we define the {\em input size of a point} 
$(v_1,\ldots,v_n)\!\in\!\Qn$ as the sum of sizes of the numerators and 
denominators of the $v_i$ (written in lowest terms), and thus extend the 
notion of input size to polynomials in $\Q[x_1,\ldots,x_n]$. Considering real 
and imaginary parts, and summing the respect sizes, we then extend the 
definition of input size further 
still to polynomials in $\Q\!\left[\sqrt{-1}\right]\![x_1,\ldots,x_n]$. \dia 
\end{dfn} 
\begin{rem} 
Note that $\size(f)$ is, up to a bounded multiple, the sum of 
the bit-sizes of all the coefficients and exponents of $f$. 
Put even more simply, assuming we write integers as usual in some fixed base, 
and we write rational numbers as fractions in lowest terms, 
$\size(f)$ is asymptotically the same as the amount of ink 
needed to write out $f$ as a sum of monomial terms. We 
extend our definition of size to a system of polynomials 
$F\!:=\!(f_1,\ldots,f_k)$ in the obvious way by setting 
$\size(F)\!:=\!\sum^k_{i=1}\size(f_i)$. Thus, for example, the 
size of an input in Theorem \ref{thm:mixed} is $\size(w)+\sum^k_{i=1}
\size(f_i)$. \dia 
\end{rem} 

\noindent 
Via a slight modification of the classical {\em Horner's Rule} \cite{cks}, it 
is easy to see that the number of ring operations needed to evaluate an 
arbitrary $f$ at an arbitrary $x\!\in\!\Cn$ easily admits an $O(\size(f)^2)$ 
upper bound.\footnote{When just counting ring operations 
we can in fact ignore the contribution of the coefficient sizes.}  
\begin{rem} 
\label{rem:speed} 
The definition of input size we use implies that our preceding 
algorithms yield a significant speed-up over earlier techniques: 
for an $n$-variate $t$-nomial $f$ of degree $d$, with $n$ and $t$ fixed, 
our algorithms have complexity polynomial in $\pmb{\log d}$. 
The best previous techniques from computational algebra, 
including recent advances on Smale's 17th Problem \cite{bps17,bcs17},  
have complexity polynomial in $\frac{(d+n)!}{d!n!}\!\geq\!\min\{d^n,n^d\}$. 
\dia 
\end{rem} 

\noindent 
{\bf Historical Notes} {\em Using convex and/or piecewise-linear 
geometry to understand solutions of algebraic equations can 
be traced back to work of Newton around 1676 \cite{newton}. 
The earliest precursor we know to the $n\!=\!1$ case of the metric estimate 
of Theorem \ref{thm:metric} can be found in work of Ostrowski from 
around 1940 \cite[Cor.\ IX, pg.\ 143]{ostrowski}.} 

{\em More recently, tropical geometry
\cite{ekl,litvinov,tropical1,bakerrumely,macsturmf} has emerged as
a rich framework for reducing deep questions in algebraic geometry to more
tractable questions in polyhedral and piecewise-linear geometry.
For instance, the combinatorial structure of amoebae was first observed 
by Gelfand, Kapranov, and Zelevinsky around 1994 \cite{gkz94}. \dia} 
\begin{rem} 
The reader may wonder why we have not considered the {\em phases} 
of the root coordinates and focussed just on norms. The phase analogue 
of an amoeba is the {\em co-amoeba}, which has only recently been studied 
\cite{herb,passare,nisse1,nisse2}. While it is known that 
the phases of the coordinates of the roots of polynomial systems 
satisfy certain equidistribution laws (see, e.g., \cite[Thm.\ 1 (pp.\ 82--83), 
Thm.\ 2 (pp.\ 87--88), and Cor.\ 3$'$ (pg.\ 88)]{few} and \cite{dgs}), there  
does not yet appear to be a phase analogue of $\trop(f)$. Nevertheless, 
we will see in Section \ref{sec:exp} that our techniques sometimes allow 
us to approximate actual complex roots, in addition to norms. \dia 
\end{rem} 

\section{Background} 
\label{sec:background} 
\subsection{Convex, Piecewise-Linear, and Tropical Geometrical Notions} 
Let us first recall the origin of the phrase ``tropical geometry'', 
according to \cite{pin}: the {\em tropical semifield} $\R_{\mathrm{trop}}$ 
is the set $\R\cup\{-\infty\}$, endowed with
with the operations $x\odot y\!:=\!x+y$ and $x\oplus y\!:=\!\max\{x,y\}$. 
The adjective ``tropical'' was coined by French computer scientists, in honor 
of Brazilian computer scientist Imre Simon, who did pioneering work 
with algebraic structures involving $\R_{\mathrm{trop}}$. 
Just as algebraic geometry relates geometric properties of zero sets 
of polynomials to the structure of ideals in commutative rings, 
tropical geometry relates the geometric properties of certain 
polyhedral complexes (see Definition \ref{dfn:complex} below) to 
the structure of ideals in $\R_{\mathrm{trop}}$. 

In our setting, we work with a particular kind of tropical variety 
that, thanks to Theorem \ref{thm:metric}, approximates $\amoeba(f)$ 
quite well. For example, one can see directly that 
$\amoeba(0)\!=\!\trop(0)\!=\!\Rn$ and, 
for any $c\!\in\!\Cs$ and $a\!\in\!\Zn$, $\amoeba(cx^a)\!=\!\trop(cx^a)\!
=\!\emptyset$. The binomial case is almost as easy. 
\begin{prop} 
\label{prop:binom} 
For any $a\!\in\!\Zn$ and non-zero complex $c_1$ and $c_2$, we have\\  
\mbox{}\hfill 
$\amoeba(c_1+c_2x^a)\!=\!\trop(c_1+c_2x^a)\!=\!
\{w\!\in\!\Rn\; | \; a\cdot w\!=\!\log|c_1/c_2|\}$.\hfill\mbox{}  
\end{prop} 

\noindent 
{\bf Proof:} If $c_1+c_2x^a\!=\!0$ then $|c_2x^a|\!=\!|c_1|$. 
We then obtain $a\cdot w\!=\!\log|c_1/c_2|$ upon taking logs and 
setting $w\!=\!\Log|x|$. 
This proves that $\amoeba(c_1+c_2x^a)$ is exactly the stated 
affine hyperplane. Similarly, since the definition of $\trop(c_1+c_2x^a)$ 
implies that we are looking for $w$ with $|c_2e^{a\cdot w}|\!=\!
|c_1|$, we see that $\trop(c_1+c_2x^a)$ defines the same hyperplane. \qed  

\smallskip 
While $\trop(f)$ and $\amoeba(f)$ are always metrically close, 
$\trop(f)$ need not be contained in, nor even have the same homotopy type as 
$\amoeba(f)$, in general.\\  
\label{ex:fancy} \addtocounter{thm}{1} 
\begin{picture}(200,170)(5,5)
\put(0,0){\epsfig{file=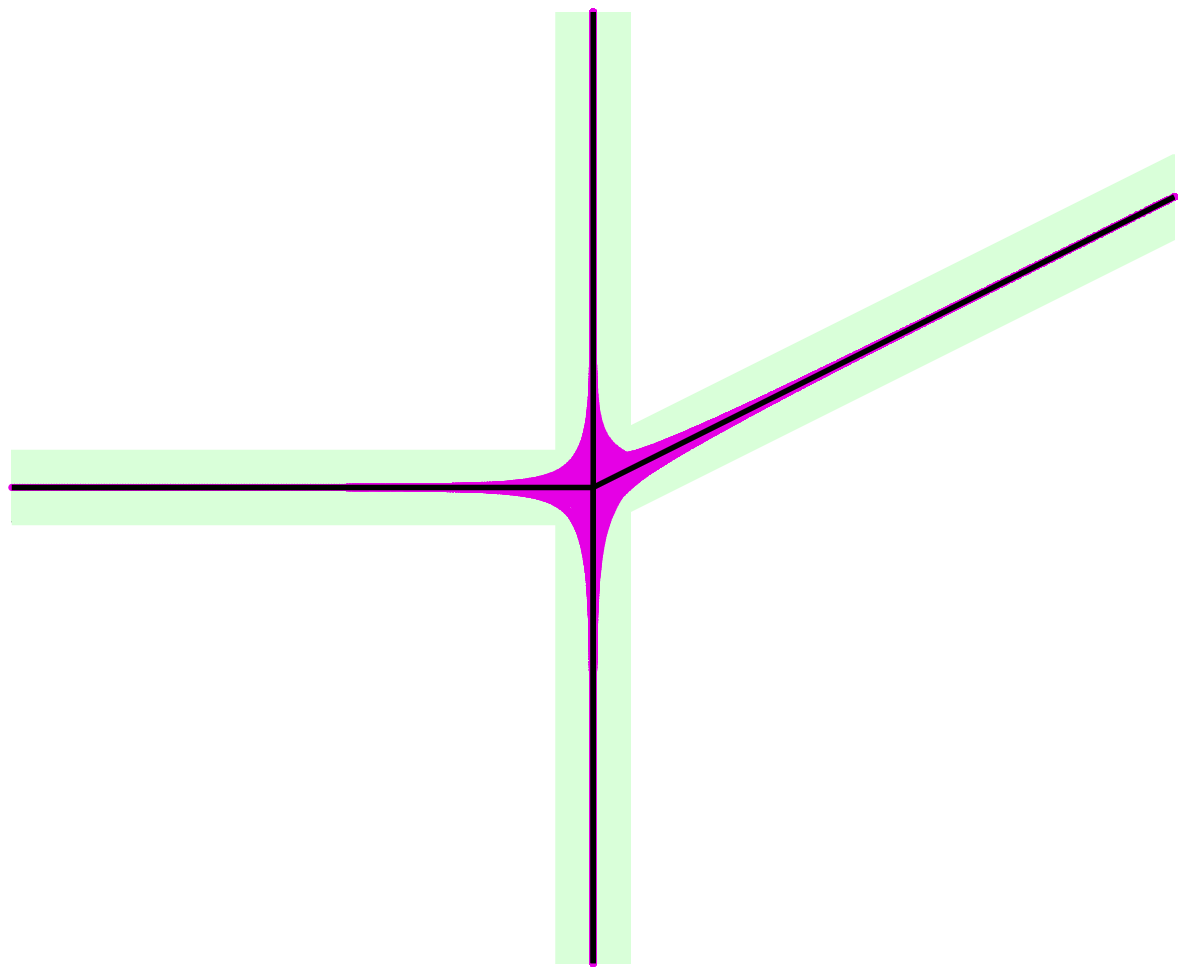,height=2.2in,clip=}\hspace{1.1cm} 
\epsfig{file=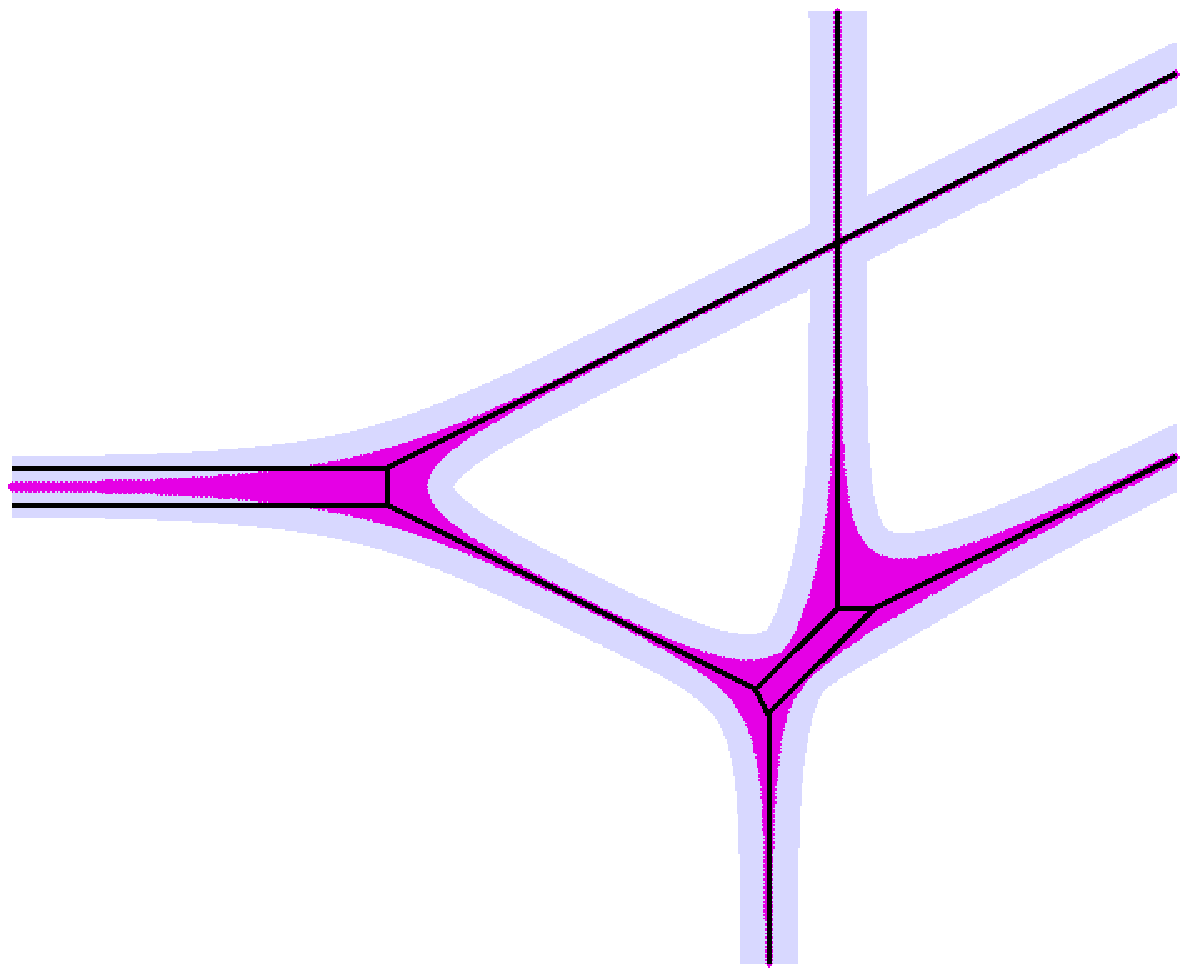,height=2.2in,clip=}}
\put(0,160){{\bf Example 2.2.}} 
\end{picture} 

\noindent 
{\em Letting $f\!:=\!1 + x^2_2 + x^4_2 + x_1 x^2_2 + x_1 x^4_2 + x^2_1 x_2 
+ x^2_1 x^2_2 + x^3_1$ and \\ 
\mbox{}\hfill 
$g\!:=\!0.1 + 0.2x^2_2 + 0.1x^4_2 + 10x_1 x^2_2 + 0.001x_1 x^4_2 + 0.01 
x^2_1 x_2 + 0.1 x^2_1 x^2_2 + 0.000005 x^3_1$\\ 
we obtain the amoebae and tropical 
varieties (and more lightly shaded neighborhoods), 
restricted to\linebreak 
$[-11,11]\times [-9,9]$, as respectively drawn on the 
left and right above. The outermost shape in the left-hand (resp.\ right-hand) 
illustration is a neighborhood of $\trop(f)$ (resp.\ $\amoeba(g)$). 

We thus see that every point of $\amoeba(f)$ (resp.\ $\trop(g)$) 
lies well within a distance of $0.65$ (resp.\ $0.49$)
of some point of $\trop(f)$ (resp.\ $\amoeba(g)$), safely within the 
distance $\log 7\!<\!1.946$ (resp.\ $13\log 7\!<\!25.3$) 
guaranteed by the second (resp.\ first) bound of Theorem \ref{thm:metric}. 
Note in particular that $\trop(g)$ has two holes while $\amoeba(g)$ has only a 
single hole.\footnote{ For our purposes, a {\em hole} of a subset 
$S\!\subseteq\!\Rn$ will simply be 
a bounded connected component of the complement $\Rn\setminus S$.} \dia}  

\smallskip 
Given any $f$, one can always easily construct a 
family of deformations whose amoebae tend to $\trop(f)$ in a suitable sense. 
This fact can be found in earlier papers of Viro and Mikhalkin, 
e.g., \cite{virologpaper,mik3}. However, employing Theorem \ref{thm:metric} 
here, we can give a $4$-line proof. 
\begin{thm}
\label{thm:deform} 
For any $n$-variate $t$-nomial $f$ written $\sum^t_{i=1} c_i x^{a_i}$, and 
$s\!>\!0$, define $f^{*s}(x)\!:=\!\sum^t_{i=1} c^s_ix^{a_i}$. Then\linebreak  
$\Delta\!\left(\frac{1}{s}\amoeba(f^{*s}),\trop(f)\right)\!\rightarrow\!0$ 
as $s\rightarrow +\infty$. 
\end{thm} 

\noindent 
\scalebox{.97}[1]{{\bf Proof:} By Theorem \ref{thm:metric}, 
$\Delta\!\left(\amoeba(f^{*s}),\trop(f^{*s})\right)\!\leq\!(2t-3)\log(t-1)$   
for all $s\!>\!0$. Since $|c_ie^{a_i\cdot w}|\!=\!|c_je^{a_j\cdot w}| 
\Longleftrightarrow$}\linebreak 
$|c_ie^{a_i\cdot w}|^s\!=\!|c_je^{a_j\cdot w}|^s$, and 
similarly when ``$=$'' is replaced by ``$>$'', we immediately obtain 
that\linebreak $\trop(f^{*s})\!=\!s\trop(f)$. So then 
$\Delta\!\left(\amoeba(f^{*s}),\trop(f^{*s})\right)\!=\!
s\Delta\!\left(\frac{1}{s}\amoeba(f^{*s}),\trop(f)\right)$   
and thus 
$\Delta\!\left(\frac{1}{s}\amoeba(f^{*s}),\trop(f)\right)\!\leq\!\frac{(2t-3)
\log(t-1)}{s}$ for all $s\!>\!0$. \qed  

\medskip 
To more easily link $\trop(f)$ with 
polyhedral geometry we will need two variations of the classical 
Newton polygon. First, let us use $\conv(S)$ to denote the 
convex hull of\footnote{i.e., smallest convex set containing...} 
a subset $S\!\subseteq\!\Rn$, $\bO\!:=\!(0,\ldots,0)$, and 
$[N]\!:=\!\{1,\ldots,N\}$. Recall also that a {\em polytope} is the 
convex hull of a finite point set, a {\em (closed) half-space} is any set of 
the form\linebreak 
$\{w\!\in\!\Rn\; | \; a\cdot w\!\leq\!b\}$ (for some $b\!\in\!\R$ and 
$a\!\in\!\Rn\setminus\{\bO\}$), and a {\em (closed) polyhedron} is any 
finite intersection of (closed) half-spaces. It is a basic fact from 
convex geometry that every polytope is a polyhedron, but not vice-versa 
\cite{grunbaum,ziegler}. 
\begin{dfn} 
\label{dfn:newt} 
Given any $n$-variate $t$-nomial $f$ written $\sum^t_{i=1} c_i x^{a_i}$, 
we define its {\em (ordinary) Newton polytope} to be
\scalebox{.95}[1]{$\newt(f)\!:=\!\conv\!\left(\{a_i\}_{i\in [t]}\right)$, 
and the {\em Archimedean Newton polytope of $f$} to be 
$\anewt(f)\!:=\!\conv\!\left(\{(a_i,-\log|c_i|)\}_{i\in [t]}\right)$.}
\linebreak   
Also, for any polyhedron $P\!\subset\!\R^N$ and $v\!\in\!\R^N$, we 
define the {\em face of $P$ with outer normal $v$} to be\linebreak  
$P^v\!:=\!\{x\!\in\!P\; | \; v\cdot x \text{ is maximized}\}$.  
The {\em dimension} of $P$, written $\dim P$, 
is simply the dimension of the smallest affine linear subspace containing 
$P$. Faces of $P$ of dimension $0$, $1$, and $\dim P-1$ are respectively 
called {\em vertices}, {\em edges}, and {\em facets}. ($P$ is called the 
{\em improper face} of $P$ and we set $\dim \emptyset\!=\!-1$.)  
Finally, we call any face of $P$ {\em lower} if and only if 
it has an outer normal $(w_1,\ldots,w_N)$ with $w_N\!<\!0$, and we 
let the {\em lower hull} of $\anewt(f)$ be the union of the lower 
faces of $\anewt(f)$. \dia 
\end{dfn} 

\noindent 
Note that $\anewt(f)$ usually has dimension $1$ greater than that of 
$\newt(f)$. $\anewt(f)$ enables us to relate $\trop(f)$ to linear 
programming, starting with the following observation.  
\begin{prop}
\label{prop:dual} 
For any $n$-variate $t$-nomial $f$, $\trop(f)$ also has the equivalent
definition\\ 
\mbox{}\hfill $\{w\!\in\!\Rn\; | \; (w,-1) \text{ is an outer normal of a 
positive-dimensional face of } \anewt(f)\}$.\hfill\mbox{} 
\end{prop}

\noindent
{\bf Proof:} The quantity $|c_ix^{a_i\cdot w}|$ being maximized at at least 
two indices $i$ is equivalent to the linear form with coefficients $(w,-1)$
being maximized at at least two difference points in
$\{(a_i,-\log|c_i|)\}_{i\in [t]}$. Since a face of a polytope is
positive-dimensional if and only if it has at least two vertices, we are done.
\qed

\begin{ex} 
\label{ex:first} 
The Newton polytope of our first example, $f\!=\!1+x^3_1+x^2_2-3x_1x_2$, 
is simply the convex hull of the exponent vectors of the monomial 
terms: $\conv(\{(0,0),(3,0),(0,2),(1,1)\})$. 
For the Archimedean Newton polytope, we 
\scalebox{.95}[1]{take the coefficients into account via an extra coordinate: 
$\anewt(f)\!=\!\conv(\{(0,0,0),(3,0,0),(0,2,0),(1,1,-\log 3)\})$.} 
\scalebox{.93}[1]{In particular, $\newt(f)$ is a triangle and $\anewt(f)$ is a 
triangular}  

\vspace{-.2cm} 
\noindent 
\begin{minipage}[b]{0.55\linewidth}
\vspace{0pt}
pyramid with base $\newt(f)\times\{0\}$ and apex 
lying beneath $\newt(f)\times\{0\}$. Note also 
that the image of the orthogonal projection of the lower hull of $\anewt(f)$ 
onto $\R^2\times \{0\}$ naturally induces a triangulation of $\newt(f)$, as 
illustrated to the right. \dia  
\end{minipage}\hspace{2.5cm}  
\begin{minipage}[b]{0.3\linewidth}
\vspace{0pt}
\epsfig{file=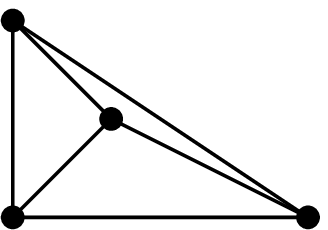,height=.7in,clip=}  
\end{minipage} 
\end{ex} 

Our last example motivates us to consider more general subdivisions and 
duality. (An outstanding reference is \cite{triang}.) Recall that a 
{\em $k$-simplex} is the convex hull of 
$k+1$ points in $\R^N$ {\em not} lying in any $(k-1)$-dimensional 
affine linear subspace of $\R^N$. A {\em simplex} is then simply a $k$-simplex 
for some $k$. 
\begin{dfn} 
\label{dfn:complex} 
A {\em polyhedral complex} is a collection of polyhedra 
$\Sigma\!=\!\{\sigma_i\}_i$ such that for all $i$ we have 
(a) every face of $\sigma_i$ is in $\Sigma$ and (b) for all $j$ we have 
that $\sigma_i\cap\sigma_j$ is a face of 
both $\sigma_i$ and $\sigma_j$. (We allow empty 
and improper faces.) The $\sigma_i$ are the {\em cells} of the complex, and 
the {\em underlying space of $\Sigma$} is 
$|\Sigma|\!:=\!\bigcup_i \sigma_i$.  

A {\em polyhedral subdivision} of a polyhedron $P$ is then simply a polyhedral 
complex $\Sigma\!=\!\{\sigma_i\}_i$ with $|\Sigma|\!=\!P$.  
We call $\Sigma$ a {\em triangulation} if and only if every $\sigma_i$ is a 
simplex. Given any finite subset $A\!\subset\!\Rn$, a 
{\em polyhedral subdivision of $A$} is then just a polyhedral subdivision 
of $\conv(A)$ where the vertices of the $\sigma_i$ all lie in $A$. 
Finally, the {\em polyhedral subdivision of $\newt(f)$ induced by 
$\anewt(f)$}, denoted $\Sigma_f$, is simply the polyhedral subdivision whose 
cells are $\{\pi(Q)\; | \; Q \text{ is a lower face of } \anewt(f)\}$, where 
$\pi : \R^{n+1}\longrightarrow \Rn$ denotes the orthogonal projection 
forgetting the last coordinate. \dia 
\end{dfn} 

Recall that a {\em (pointed polyhedral) cone} is just the set of all 
nonnegative linear combinations of a finite set of points. Such cones 
are easily seen to always be polyhedra \cite{grunbaum,ziegler}. Recall also 
that a {\em bijection}, $\phi$, between two finite sets $A$ and $B$ is just a  
function $\phi : A\longrightarrow B$ such that the cardinalities of 
$A$, $B$, and $f(A)$ are all equal. 
\begin{ex} 
The illustration from Example \ref{ex:first} shows a triangulation 
of the point set $\{(0,0),(3,0),(0,2),(1,1)\}$ which happens to be 
$\Sigma_f$ for $f\!=\!1+x^3_1+x^2_2-3x_1x_2$. More to the point, it is 
easily checked that the outer normals to a face of dimension $k$ of $\anewt(f)$ 
form a cone of dimension $3-k$. In this way, thanks to the 
natural partial ordering of cells in any polyhedral complex by inclusion, 
we get an order-reversing bijection between the cells of $\Sigma_f$ and 
pieces of $\trop(f)$. \dia 
\end{ex} 
That $\trop(f)$ is always a polyhedral complex follows directly  
from Proposition \ref{prop:dual} above. It is then easy to show that 
there is always an order-reversing bijection between the  
cells $\Sigma_f$ and the cells of $\trop(f)$ --- an incarnation of 
{\em polyhedral duality} \cite{ziegler}. 
\begin{ex} 
We illustrate the preceding order-reversing bijection of cells through 
our first three tropical varieties, and corresponding subdivisions $\Sigma_f$ 
of $\newt(f)$, below: \\ 
\epsfig{file=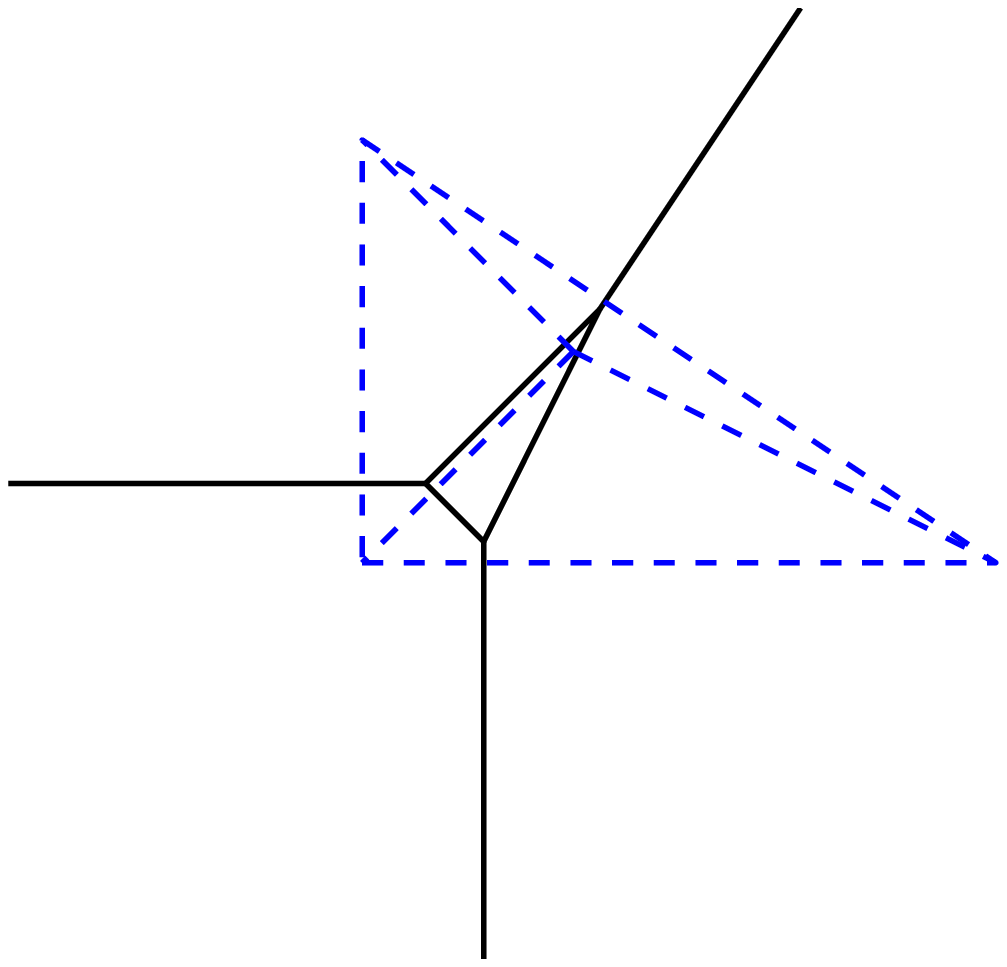,height=2in,clip=}
\epsfig{file=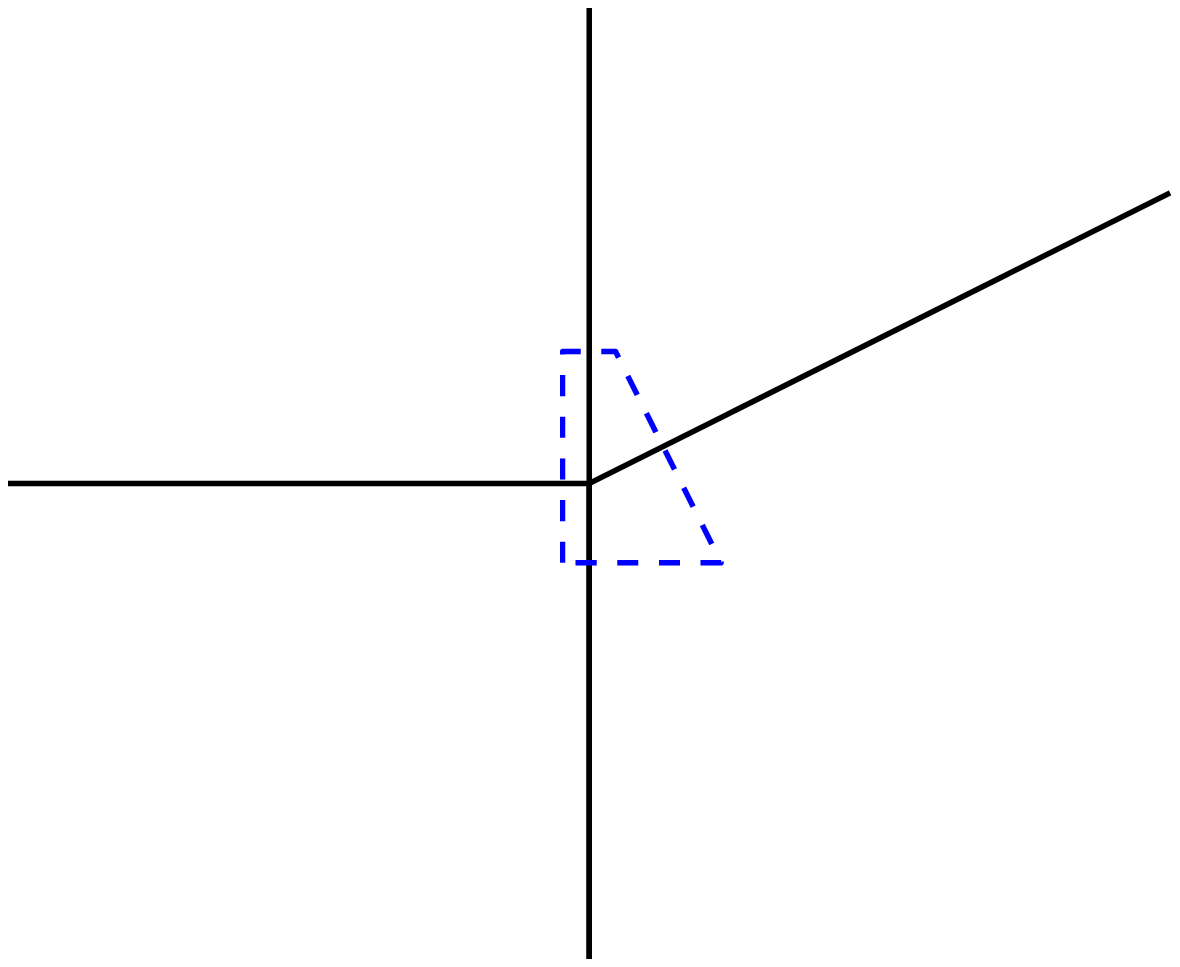,height=2in,clip=}
\epsfig{file=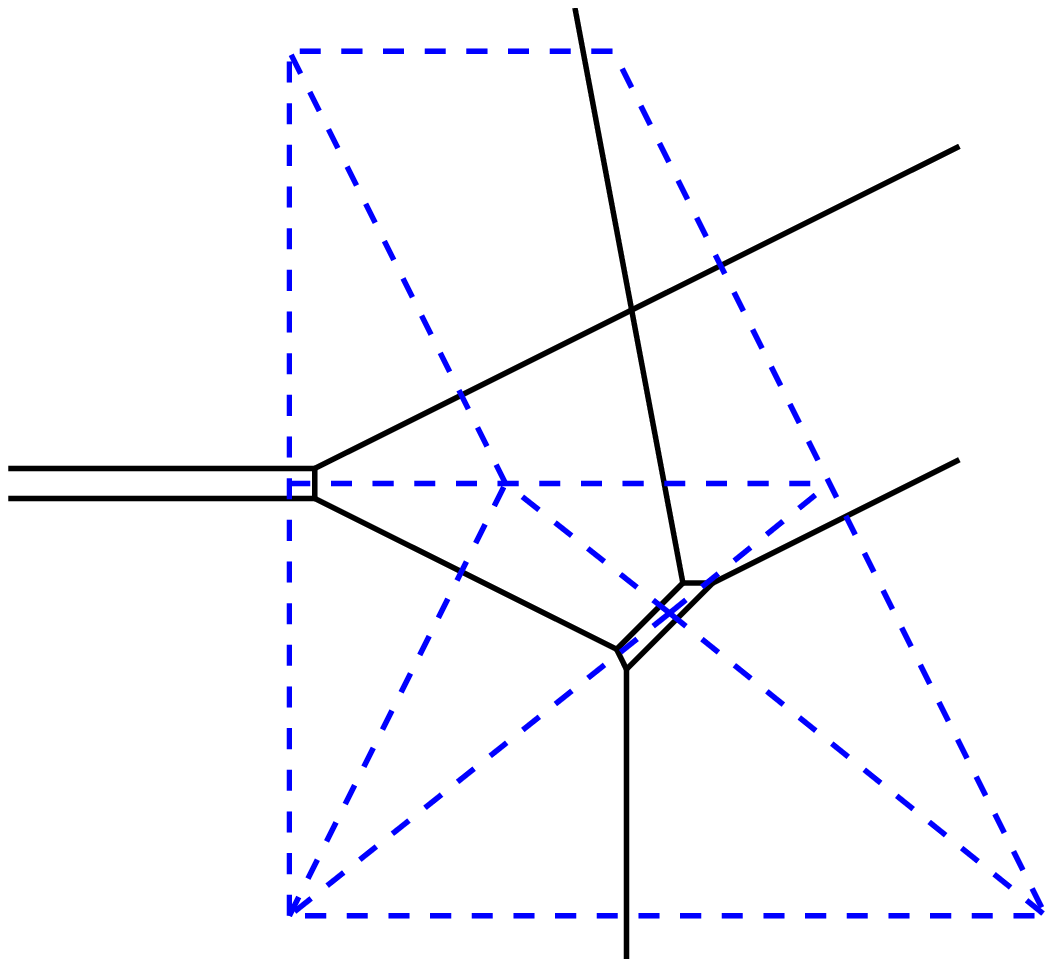,height=2in,clip=}\\ 
Note that the vertices of $\trop(f)$ correspond bijectively to the 
$2$-dimensional cells of $\Sigma_f$, and the $1$-dimensional cells of 
$\trop(f)$ correspond 
bijectively to the edges of $\Sigma_f$. (In particular, the rays of 
$\trop(f)$ are perpendicular to the edges of $\newt(f)$.) 
Note also that the vertices of 
$\Sigma_f$ correspond bijectively to connected components of the {\em 
complement} $\R^2\!\setminus\!\trop(f)$. We have taken the liberty of 
slightly distorting the right-most illustration to make the bijections 
clearer. \dia 
\end{ex} 

\subsection{The Complexity of Linear Programming} 
\label{sub:lo} 
Let us first point out that \cite{papa,cxity,sipser} are outstanding 
references for further background on the classical Turing model and 
$\np$-completeness. Let us now focus on some well-known late-20$\thth$ 
century results on the complexity of linear optimization. These 
results are covered at much greater length in \cite{schrijver,grotschel}. 
\begin{dfn}  
Let $\R^N_\geq\!:=\!\{(x_1,\ldots,x_N)\!\in\!\Rn\; | \; x_1,\ldots,
x_N\!\geq\!0\}$ denote the nonnegative orthant. 
Given $M\!\in\!\R^{k\times N}$\linebreak 
with linearly independent rows, $c\!=\!(c_1,\ldots,c_N)\!\in\!\R^N$, and 
$b\!=\!(b_1,\ldots,b_k)\!\in\!\R^k$, the {\em (standard form) linear 
optimization problem} $\cL(M,b,c)$ is the following problem:\\ \mbox{}\\  
\mbox{}\hspace{2in}Maximize $c\cdot x$ subject to:\\ 
\mbox{}\hspace{2in} $Mx=b$ \\ 
\mbox{}\hspace{2in} $x\!\in\!\R^N_\geq$ 

\smallskip 
\noindent   
We then define $\size(\cL(M,b,c))\!:=\!\size(M)+\size(b)+\size(c)$. 
The set of all $x\!\in\!\R^N_\geq$ satisfying $Mx\!=\!b$ is the 
{\em feasible region} of $\cL(M,b,c)$. We call $\cL(M,b,c)$ 
{\em infeasible} if and only if there is {\em no} $x\!\in\!\R^N_\geq$ 
satisfying $Mx\!=\!b$. Finally, if $\cL(M,b,c)$ is feasible but does not 
admit a well-defined maximum, then we call $\cL(M,b,c)$ {\em unbounded}. \dia 
\end{dfn} 
\begin{thm} 
\label{thm:lp} 
Given any linear optimization problem $\cL(M,b,c)$ as defined above, 
we can decide infeasibility,\linebreak 
\scalebox{.96}[1]{unboundedness, or (if $\cL(M,b,c)$ is feasible) find an 
optimal solution $x^*$, all within time polynomial in 
$\size(\cL(M,b,c))$.}\linebreak 
\scalebox{1}[1]{In particular, if $\cL(M,b,c)$ is feasible, we can find an 
optimal solution $x^*$ of size polynomial in $\size(\cL(M,b,c))$. \qed}  
\end{thm} 

\noindent 
Theorem \ref{thm:lp} goes back to work of Khachiyan in the late 1970s on the 
{\em Ellipsoid Method}, building upon earlier work of Shor, Yudin, and 
Nemirovskii \cite{schrijver}. Since then, {\em Interior Point Methods} have 
emerged as one of the most practical methods attaining the complexity 
bound asserted in Theorem \ref{thm:lp}. For simplicity, we will not 
focus on the best current complexity bounds, since we 
simply want to prove polynomiality for our algorithms in this paper. 
Further discussion on improved complexity bounds for linear optimization can be 
found in \cite{monteiro}. 

Any system of linear inequalities, at the expense of a minor 
increase in size, is essentially equivalent to the feasible 
region of some $\cL(M,b,c)$. In what follows, $Mx\!\leq\!b$ is 
understood to mean that $M_1\cdot x\!\leq\!b_1,\ldots,M_k\cdot x\!\leq\!b_k$ 
all hold, where $M_i$ denotes the $i\thth$ row of $M$. 
\begin{prop} 
\label{prop:convert} 
Given $M\!\in\!\R^{k\times N}$ and any collection of inequalities of the 
form $Mx\!\leq\!b$, there is a standard form linear optimization problem 
$\cL(\bar{M},\bar{b},\bO)$, satisfying 
$\size(\cL(\bar{M},\bar{b},\bO))\!\leq\!2(\size(M)+\size(b))+k$, 
that is feasible if and only if $\{x\!\in\!\Rn\; | \; Mx\!\leq\!b\}$ is 
non-empty. \qed 
\end{prop} 
There is thus no loss of generality in restricting to standard form. 

We will frequently work with polyhedra given explicitly in the form 
$P\!=\!\{x\!\in\!\Rn\; | \; Mx\!\leq\!b\}$ (usually called {\em 
$\cH$-polytopes}), and use Proposition 
\ref{prop:convert} and Theorem \ref{thm:lp} together to rapidly decide various 
basic questions about $P$. For instance, 
we call a constraint $M_i\cdot x\!\leq\! b_i$ of $Mx\!\leq\!b$ 
{\em redundant} if and only if 
the corresponding row of $M$ can be deleted from $M$ without affecting $P$. 
\begin{lem} 
\label{lemma:red} 
Given any system of linear inequalities $Mx\!\leq\!b$ we can, 
in time polynomial in\linebreak $\size(M)+\size(b)+\size(c)$,  
find a submatrix $M'$ of $M$ (and a subvector $b'$ obtained 
by deleting the corresponding\linebreak entries from $b$) such that 
$\{x\!\in\!\R^N\; | \; M'x\!\leq\!b'\}\!=\!\{x\!\in\!\R^N\; | \; 
M'x\!\leq\!b'\}$ and $M'x\!\leq\!b'$ has no redundant 
constraints. \qed 
\end{lem} 

\noindent 
The new set of inequalities $M'x\!\leq\!b'$ is called an {\em irredundant 
representation} of $Mx\!\leq\!b$. 

A deep subtlety underlying linear optimization is whether 
$\cL(M,b,c)$ can be solved in {\em strongly} polynomial-time, i.e., 
is there an analogue of Theorem \ref{thm:lp} where we instead 
count {\em arithmetic operations} to measure complexity, and obtain 
complexity polynomial in $k+N$? 

One of the first successful algorithms for linear optimization --- the {\em 
Simplex Method} --- has {\em arithmetic} complexity $O(N^k)$, and there are 
now variations of the Simplex Method (using sophisticated {\em pivoting rules}) 
that attain arithmetic complexity sub-exponential in $k$. (It was also 
discovered in the 1970s by Borgwardt and Smale that the simplex method is 
strongly polynomial {\em provided one averages over a suitable distribution of 
inputs} \cite{schrijver}.) Strong polynomiality remains an important open 
problem and is in fact Problem 9 on Fields Medalist Steve Smale's list of 
mathematical problems for the 21$\stst$ Century \cite{21a,21b}.  

These issues are actually relevant to polynomial system solving since 
the linear optimization problems we ultimately solve will have 
{\em irrational ``right-hand sides''}: $b$ will usually be a 
(rational) linear combination of logarithms of integers in our setting. 

In particular, as is well-known in Diophantine Approximation \cite{baker}, 
it is far from trivial to efficiently decide the sign of such an 
irrational number. This problem is also easily seen to be equivalent to 
deciding inequalities of the form $\alpha^{\beta_1}_1\cdots 
\alpha^{\beta_N}_N\!\stackrel{?}{>}\!1$, where the $\alpha_i$ and 
$\beta_i$ are integers. Note, in particular, that while the number of 
arithmetic operations necessary to decide such an inequality is 
easily seen to be $O((\sum^N_{i=1}\log|\beta_i|)^2)$ (via the classical 
binary method of exponentiation), taking bit-operations into account 
naively results in a problem that appears to have complexity {\em exponential}  
in $\log|\beta_1|+\cdots+\log|\beta_N|$. Fortunately, another Fields Medalist, 
Alan Baker, made major progress on this problem in the late 20th century. 

\subsection{Irrational Linear Optimization and Approximating Logarithms Well 
Enough} 
\label{sub:log} 
Recall the following result on comparing monomials in rational numbers.
\begin{thm} \cite[Sec.\ 2.4]{brs}
\label{thm:linlog}
Suppose $\alpha_1,\ldots,\alpha_n\!\in\!\Q$ are positive and
$\beta_1,\ldots,\beta_n\!\in\!\Z$. Also let $A$ be the
maximum of the numerators and denominators of the $\alpha_i$ (when
written in lowest terms) and $B\!:=\!\max_i\{|\beta_i|\}$. Then, within\\
\mbox{}\hfill $O\!\left(n30^n\log(B)(\log \log B)^2\log\log\log(B)
(\log(A)(\log \log A)^2\log\log\log A)^n\right)$\hfill \mbox{} \\
bit operations, we can determine the sign of $\alpha^{\beta_1}_1\cdots 
\alpha^{\beta_n}_n-1$. \qed
\end{thm}

\noindent
While the underlying algorithm is a simple application of
Arithmetic-Geometric Mean Iteration (see, e.g.,
\cite{dan}), its complexity bound hinges on a deep estimate
of Nesterenko \cite{nesterenko}, which in turn refines seminal
work of Matveev \cite{matveev} and Alan Baker \cite{baker} on
linear forms in logarithms. Whether the dependence on $n$ in the 
bound above can be improved to polynomial is a very deep open 
question related to the famous $abc$-Conjecture \cite{bakerabc,nitaj}.  

Via the Simplex Method, or even a brute force search through all 
basic feasible solutions of $\cL(M',b',c')$, we can obtain the 
following consequence of Theorems \ref{thm:lp} and \ref{thm:linlog}.   
\begin{cor} 
\label{cor:irrat} 
Suppose $n$ is fixed, $k\!\leq\!n$, $M\!\in\!\Q^{k\times n}$, 
and $b_i\!:=\!\log|\beta_i|$ with 
$\beta_i\!\in\!\Q^*$ for all $i\!\in\![k]$, and we set\linebreak  
$b\!:=\!(b_1,\ldots,b_k)$. Then we can decide feasibility  
for $Mx\!\leq\!b$, and compute an irredundant representation\linebreak  
$M'x\!\leq\!b'$ for $Mx\!\leq\!b$, in time polynomial in $\size(M)+\size(b)$. 
\qed 
\end{cor} 

\noindent 
The key trick behind the proof of Corollary \ref{cor:irrat} is that, 
after converting to standard form, any basic feasible solution of 
the underlying linear optimization problem has all its irrationalities 
concentrated on the right-hand side. In particular, standard linear 
algebra bounds tell us that the right-hand side involves a linear 
combination of logarithms with coefficients of size polynomial in the 
input size. 

\section{Tropical Start-Points for Numerical Iteration and an Example} 
\label{sec:exp} 
We begin by outlining a method for picking start-points for 
Newton Iteration (see, e.g., \cite[Ch.\ 8]{bcss} for a modern 
perspective) and Homotopy Continuation \cite{huangli,sw,verschelde,leeli,bhsw}. 
While we do not discuss 
these methods for solving polynomial equations in detail, let us point 
out that Homotopy Continuation (combined with Smale's $\alpha$-Theory 
for certifying roots \cite{bcss,bhsw}) is currently the fastest and most 
reliable method for numerically solving polynomial systems in 
complete generality. Other important methods include Resultants \cite{emican} 
and Gr\"obner Bases \cite{faugere}. However, while these 
alternative methods are of great importance in certain algebraic and 
theoretical applications \cite{resheight,faugerecrypto}, Homotopy Continuation 
is currently the method of choice for practical large-scale numerical 
computation.

While the \fbox{boxed} steps below admit a simple and easily 
parallelizable brute-force search, they form the 
portion of the algorithm that is the most challenging to speed up to complexity 
polynomial in $n$. 
\begin{algor}
\label{algor:start}
(Coarse Approximation to Roots with Log-Norm Vector Near a Query Point)\\   
{\sc Input.} Polynomials $f_1,\ldots,f_n\!\in\!\C\!\left[x^{\pm 1}_1,
\ldots,x^{\pm 1}_n\right]$, with 
$f_i(x)\!=\!\sum^{t_i}_{j=1} c_{i,j} x^{a_j(i)}$ 
a $t_i$-nomial for all $i$, and a query point\\ 
\mbox{}\hspace{1.1cm}$w\!\in\!\Rn$.\\   
{\sc Output.} An ordered $n$-tuple of sets of indices $(J_i)^n_{i=1}$ such 
that $g_i\!:=\!\sum_{j\in J_i}c_{i,j}x^{a_j(i)}$ 
is a sub-summand of $f_i$, and the\\ 
\mbox{}\hspace{1.3cm}roots of $G\!:=\!(g_1,\ldots,g_n)$ 
are near the roots of $F\!:=\!(f_1,\ldots,f_n)$ with log-norm vector near 
$w$.\\   
{\sc Description.}\\  
1. Let $\sigma_w$ be the closure of the unique cell of 
$\Rn\setminus\bigcup^n_{i=1} \trop(f_i)$ or 
$\bigcap\limits_{\trop(f_i)\ni w}\trop(f_i)$ containing $w$. \\ 
\fbox{2.} If $\sigma_w$ has no vertices in $\bigcap^n_{i=1}\trop(f_i)$ then 
output an irredundant collection of facet inequalities for $\sigma_w$,\\ 
\mbox{}\hspace{.6cm}output {\tt ``There are no roots of $F$ in $\sigma_w$.''}, 
and {\tt STOP}.\\ 
\fbox{3.} Otherwise, fix a vertex $v$ of $\sigma_w$ and, for each 
$i\!\in\![n]$, let $E_i$ be any edge of $\anewt(f)$ 
generating a facet of\\ 
\mbox{}\hspace{.6cm}$\trop(f_i)$ containing $v$. \\ 
4. For all $i\!\in\![n]$, let $J_i\!:=\!\{j\; | \; (a_j(i),
-\log|c_{i,j}|)\!\in\!E_i\}$.\\ 
5. Output $(J_i)^n_{i=1}$. \qed 
\end{algor} 
\begin{rem} 
The output system $G$ is useful because, with probability $1$ (for most 
reasonable distributions on the coefficients), all the $g_i$ are 
{\em binomials}, and binomial systems are particularly easy to solve: 
they are equivalent to linear equations in the logarithms of the original 
variables. In particular, an $n\times n$ binomial system output 
by our algorithm always results in a collection of roots all sharing a 
{\em single} vector of norms. 

The connection to Newton Iteration is then easy to state: use any root of 
$G$ as a start-point $z(0)$ for the iteration\\ 
\mbox{}\hfill $z(n+1)\!:=\!
z(n)-\mathrm{Jac}(F)^{-1}|_{z(n)}F(z(n))$.\hfill\mbox{}\\ 
The connection to Homotopy 
Continuation is also simple: use the pair $(G,\zeta)$ (for any 
root $\zeta$ of $G$) to start a path converging (under the usual 
numerical conditioning assumptions on whatever predictor-corrector method 
one is using) to a root of $F$ with log-norm vector near $w$. Note that 
it is safer to do the extra work of Homotopy Continuation, but there 
will be cases where the tropical start-points from 
Algorithm \ref{algor:start} are sufficiently good that Newton Iteration is 
enough to converge to a true root. 

Note in particular that we have the freedom to follow as few start-points, 
or as few paths, as we want. When our start-points (resp.\ paths) indeed 
converge to nearby roots, we obtain a tremendous savings over having to follow 
{\em all} start-points (resp.\ paths). \dia 
\end{rem} 
\begin{dfn} 
\label{dfn:mix} 
Following the notation of Theorem \ref{thm:mixed} and Algorithm 
\ref{algor:start}, we call a vertex $v$ of $\sigma_w$ {\em mixed} 
if and only if it lies in $\bigcap^n_{i=1}\trop(f_i)$. \dia 
\end{dfn} 
\begin{ex} 
Let us make a $2\times 2$ polynomial system out of our first and 
third examples:\\ 
\mbox{}\hfill $f_1\!:=\!1+x^3_1+x^2_2-3x_1x_2$\hfill \mbox{}\\ 
\mbox{}\hfill $f_2\!:=\!0.1 + 0.2x^2_2 + 0.1x^4_2 + 10x_1 x^2_2 
+ 0.001x_1 x^4_2 + 0.01 
x^2_1 x_2 + 0.1 x^2_1 x^2_2 + 0.000005 x^3_1$  \hfill \mbox{}\\ 
\mbox{}\hfill \epsfig{file=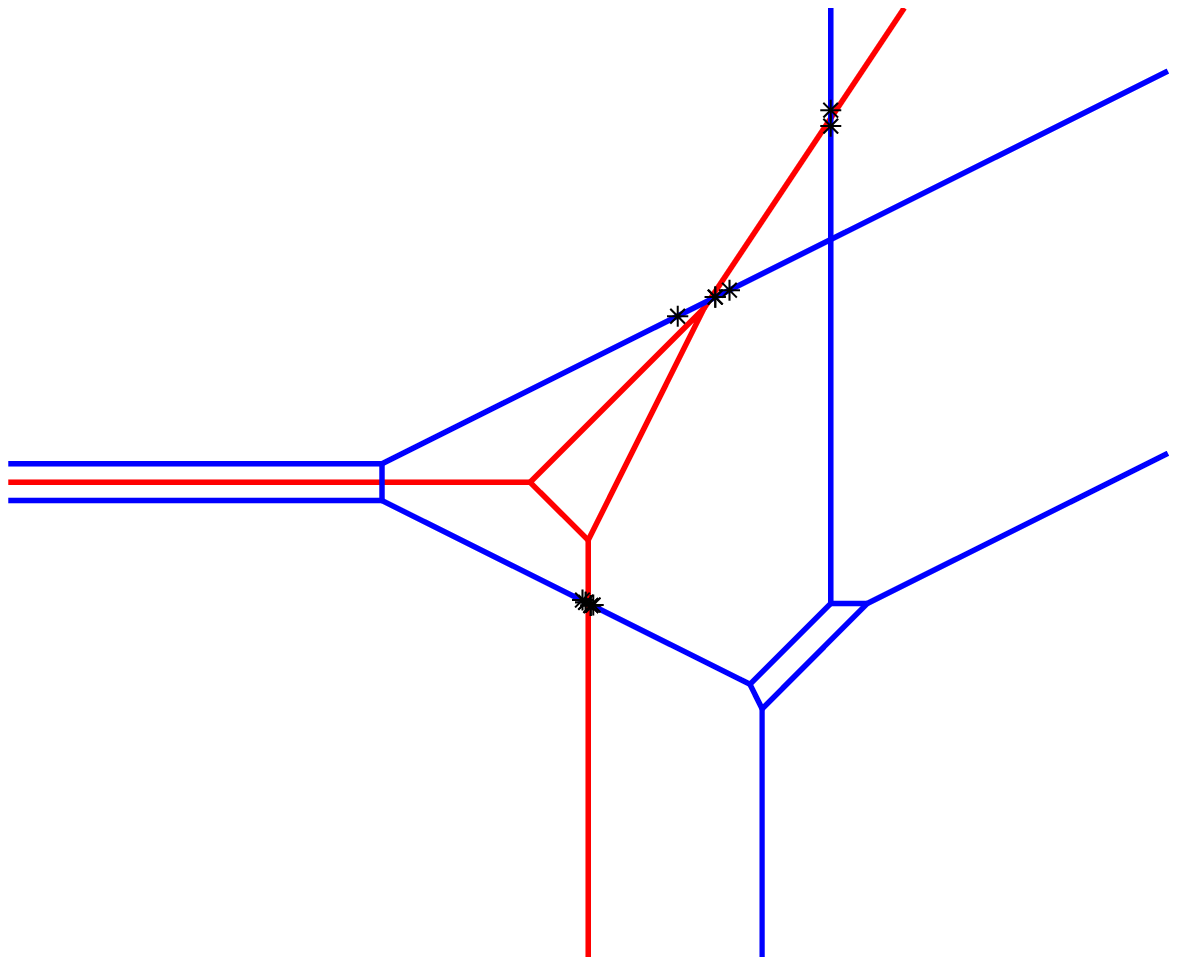,height=2.2in,clip=}\hspace{1in} 
\epsfig{file=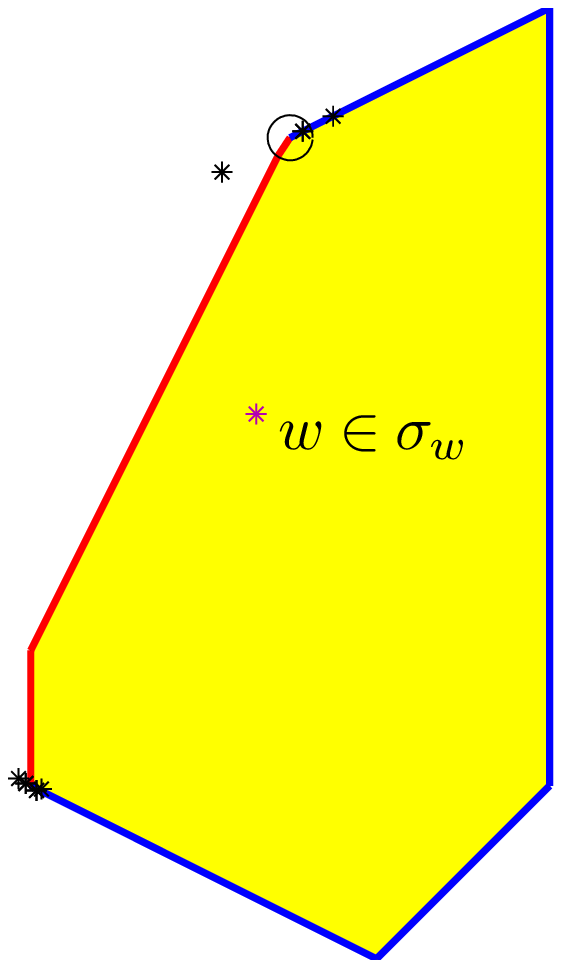,height=2.2in,clip=}
\hfill\mbox{}\\ 
The system $F\!:=\!(f_1,f_2)$ has exactly $12$ roots in $(\Cs)^2$, 
the coordinate-wise log-norms of which form the small clusters near certain 
intersections of $\trop(f_1)$ and $\trop(f_2)$.\footnote{The root 
count was verified via an exact Gr\"obner basis calculation using 
the commercial software package {\tt Maple 14} . Numerical 
approximation of the log-norm vectors to accuracy $10^{-4}$ per 
coordinate was then done via the 
publically available package {\tt Bertini} \cite{bhsw}, using 
default settings. Both calculations took a fraction of a second. 
The cell $\sigma_w$ was computed via {\tt Matlab 7.11.0 (R2010b)}. } In 
particular, there is a heptagonal cell, which we have magnified, with 
$2$ vertices close to the log-norm vectors of some of the roots.
This cell, which looks hexagonal because it has 
a pair of vertices that are too close to distinguish visually, happens to be 
$\sigma_w$ for $w\!=\!(2,1)$. Note that $\sigma_w$ has exactly $2$ {\em mixed} 
vertices. 

Applying Algorithm \ref{algor:start} to our $(f_1,f_2,w)$ we then have 
$2$ possible outputs, depending on which mixed vertex of $\sigma_w$ we 
pick. The output corresponding to the circled vertex is the 
pair of index sets $(\{2,3\},\{3,4\})$. More concretely, 
Algorithm \ref{algor:start} alleges that the system \\ 
\mbox{}\hfill $G\!:=\!(g_1,g_2)\!:=\!(x^3_1+x^2_2,0.1x^4_2+10x_1x^2_2)$ 
\hfill \mbox{}\\ 
has roots with log-norm vector near a log-norm vector of a root of $F$ that 
is in turn close to $w$. Indeed, 
the sole log-norm vector coming from the roots of $G$ is $\left(\log 10,
\frac{3}{2}\log 10\right)$ 
and the roots themselves are $\{(\pm 10,\sqrt{\mp 1000})\}$ (with 
both values of the square root allowed). All $4$ roots 
in fact converge (under Newton iteration, with no need for Homotopy 
Continuation) to true roots of $F$: $(-10,1000^{1/2})$ and 
$(-10,-1000^{1/2})$ respectively converge to the roots of $F$ 
with closest and third closest log-norm vector to $w$. 
The other two roots of $G$ converge to a conjugate pair of 
roots of $F$ with log-norm vector $(2.4139,3.5103)$ (to $4$ decimal places) 
lying in the small circle.  
\dia 
\end{ex} 

\begin{rem}
The cell $\sigma_w$ from Step 1 can be found in polynomial-time, thanks to
Theorem \ref{thm:region}, and its underlying algorithm contained in
Corollary \ref{cor:irrat}. 

As for Steps 2 and 3, thanks to duality, the facets of $\trop(f_i)$ correspond
exactly to lower edges of $\anewt(f_i)$. So, to find the vertex
$v$ (or decide that it doesn't exist), it suffices to do a brute-force
search through all $n$-tuples of lower edges, one coming from each of
$\anewt(f_1),\ldots,\anewt(f_n)$. This particular kind of geometric
computation has its origins in the algorithmic study of {\em mixed volume} 
\cite{emican,leeli}. 
There are various ways of speeding up this search and there is much
interesting computational geometry to be studied in this direction. \dia
\end{rem}

Let us be clear that we have not yet proved a metric guarantee for Algorithm 
\ref{algor:start} in the spirit of Theorem \ref{thm:metric}. Rigorous results 
in this direction, as well as a broad experimental understanding of 
our techniques, are of the utmost importance and we hope to address these  
points in the near future. 
\begin{rem} 
We have intentionally written Algorithm \ref{algor:start} in terms of a 
more general class of inputs than necessary for our examples. For such 
general inputs, it makes more sense to measure complexity in terms of 
arithmetic operations instead of bit operations. \dia 
\end{rem} 
\begin{rem} 
The reader should be aware that while we have relied upon Diophantine 
approximation and subtle aspects of the Simplex Method to prove our 
bit-complexity bounds, one can certainly be more flexible when using our 
approach in practical, floating-point computations. For instance, 
heuristically, it appears that one can get away with less accuracy than 
stipulated by Theorem \ref{thm:linlog} when comparing linear 
combinations of logarithms. Similarly, one should feel free to 
use the fastest (but still reasonably accurate) algorithms for 
linear optimization when applying our methods to large-scale 
polynomial systems. \dia 
\end{rem} 

\section{Proof of Theorem \ref{thm:region}}
\label{sec:region} 
Using $t-1$ comparisons, we can isolate all indices $i$ such that 
$\max_i|c_ie^{a_i\cdot w}|$ is attained. Thanks to Theorem \ref{thm:linlog}, 
this can be done in polynomial-time. We then obtain, say, $J$ equations of the 
form $a_i\cdot w\!=\!-\log|c_i|$ and $K$ inequalities of the 
form $a_i\cdot w\!>\!-\log|c_i|$ or $a_i\cdot w\!<\!-\log|c_i|$. 

Thanks to Lemma \ref{lemma:red}, combined with Corollary \ref{cor:irrat}, 
we can determine the exact cell of $\trop(f)$ containing $w$ if $J\!\geq\!2$. 
Otherwise, we obtain the unique cell of $\Rn\!\setminus\!\trop(f)$ containing 
$w$. Note also that an $(n-1)$-dimensional face of either kind of 
cell must be the dual of an edge of $\anewt(f)$. Since every edge 
has exactly $2$ vertices, there are at most $t(t-1)/2$ such $(n-1)$-dimensional 
faces, and thus $\sigma_w$ is the intersection of at most $t(t-1)/2$ 
half-spaces. So we are done. \qed  
\begin{rem}
Theorem \ref{thm:region} also generalizes an earlier complexity
bound for deciding membership in $\trop(f)$ from \cite{aknr}. \dia
\end{rem}

\section{Proof of Theorem \ref{thm:dist}}
\label{sec:dist} 
Since $\trop(f)$ and $\amoeba(f)$ are closed, $\Delta(w,\trop(f))\!=\!
|w-v|$ for some point $v\!\in\!\trop(f)$ and 
$\Delta(w,\amoeba(f))\!=\!|w-u|$ for some point $u\!\in\!\amoeba(f)$.  

Now, by the second upper bound of Theorem \ref{thm:metric}, 
there is a point $v'\!\in\!\trop(f)$ within distance $\log(t-1)$ 
of $u$. Clearly, $|w-v|\!\leq\!|w-v'|$. Also, by the Triangle Inequality, 
$|w-v'|\!\leq\!|w-u|+|u-v'|$. So then,\\ 
\mbox{}\hfill $\Delta(w,\trop(f))\!\leq\!
\Delta(w,\amoeba(f))+\log(t-1)$,\hfill\mbox{}\\
and thus $\Delta(w,\amoeba(f))-\Delta(w,\trop(f))\!\geq\!-\log(t-1)$.

Similarly, by the first upper bound of Theorem \ref{thm:metric},
there is a point $u'\!\in\!\amoeba(f)$ within distance\linebreak 
$(2t-3)\log(t-1)$ 
of $v$. Clearly, $|w-u|\!\leq\!|w-u'|$. Also, by the Triangle Inequality, 
$|w-u'|\!\leq\!|w-v|+|v-u'|$. So then, $\Delta(w,\amoeba(f))\!\leq\!
\Delta(w,\trop(f))+(2t-3)\log(t-1)$, and thus\\ 
\mbox{}\hfill $\Delta(w,\amoeba(f))
-\Delta(w,\trop(f))\!\leq\!(2t-3)\log(t-1)$.\hfill\mbox{}\\ 
So our first assertion is proved.  

Now, if $f$ has coefficients with rational real and imaginary parts, 
Theorem \ref{thm:region} tells us that we have an explicit 
description of $\sigma_w$ as the intersection of a 
number of half-spaces polynomial in the input size. 
Moreover, the bit-sizes of the coefficients of the underlying inequalities 
are also polynomial in the input size.  So we can compute the distance $D$ 
from $w$ to $\trop(f)$ by finding which 
facet of $\sigma_w$ has minimal distance to $w$. The 
distance from $w$ to any such facet can be computed in polynomial-time via the 
classical formula for distance between a point and an affine hyperplane, and 
Theorem \ref{thm:linlog}: \\ 
\mbox{}\hfill $\Delta(w,\{x\; | \; \alpha\cdot x=\beta\}) = 
\frac{|\alpha\cdot w|-\sign(\alpha\cdot w)\beta}{|\alpha|}$ \hfill \mbox{}\\  
In particular, we may efficiently approximate $D$ by
efficiently approximating the underlying square-roots and logarithms.
The latter can be accomplished by Arithmetic-Geometric Iteration, as detailed 
in \cite{dan}. So our statement on leading bits is proved.  

The final assertion then follows easily: we merely decide 
whether $\Delta(w,\trop(f))$ strictly exceeds $\log(t-1)$ or not, 
via the algorithm we just outlined. Thanks to our initial observations using 
the Triangle Inequality, it is clear that Output (b) 
or Output (a) occurs according as $\Delta(w,\trop(f))\!>\!\log(t-1)$ or 
not. \qed 

\section{Proving of Theorem \ref{thm:mixed}}
\label{sec:mixed} 
\subsection{Fast Cell Computation: Proof of the First Assertion} 
First, we apply Theorem \ref{thm:region} to $(f_i,w)$ 
for each $i\!\in\![k]$ to find which $\trop(f_i)$ 
contain $w$. 

If $w$ lies in no $\trop(f_i)$, then we simply use Corollary \ref{cor:irrat} 
(as in our proof of Theorem \ref{thm:region}) to find an explicit description 
of the closure of the cell of $\Rn\!\setminus\!
\bigcup^k_{i=1}\trop(f_i)$ containing $w$. Otherwise, 
we find the cells of $\trop(f_i)$ (over those $i$ 
with $\trop(f_i)$ containing $w$) that contain $w$. 
Then, applying Corollary \ref{cor:irrat} once again, we 
find the unique cell of $\bigcap\limits_{\trop(f_i)\ni w}\trop(f_i)$ 
containing $w$. 

Assume that $f_i$ has exactly $t_i$ monomial terms for all $i$. In either of 
the preceding cases, the total number of half-spaces 
involved is no more than $\sum^k_{i=1}t_i(t_i-1)/2$. So the over-all 
complexity of our redundancy computations is polynomial in the 
input size and we are done. \qed 

\subsection{Hardness of Detecting Mixed Vertices: Proving the 
Second Assertion} 
It will clarify matters if we consider a related $\np$-hard 
problem for rational polytopes first, before moving on to 
cells with irrationalities. 
\subsubsection{Preparation over $\Q$} 
In the notation of Definition \ref{dfn:mix}, let us first consider the 
following decision problem. We assume all polyhedra are given explicitly as 
finite collections of rational linear inequalities, with size defined as in 
Section \ref{sub:lo}. 

\smallskip 
\noindent 
{\sc Mixed-Vertex}:\\
{\em Given $n\in \N$ and polyhedra $P_1,\ldots,P_n$ in $\R^n$, 
does $P:=\bigcap_{i=1}^n P_i$ have a mixed vertex? \qed} 

\smallskip 
While {\sc Mixed-Vertex} can be solved in polynomial time when the dimension is 
fixed, we will show that, for $n$ varying, the problem is $\np$-complete, 
even when restricting to the case where all polytopes are full-dimensional and 
$P_1,\ldots,P_{n-1}$ are axes-parallel bricks.

Let $e_i$ denote the $i\thth$ standard basis vector in $\R^n$.
Also, given $\alpha\!\in\!\Rn$ and $\beta\!\in\!\R$, we will use the following 
notation for certain {\em hyperplanes} and {\em halfspaces} in $\R^n$ 
determined by $\alpha$ and $\beta$: 
$$
H_{(\alpha,\beta)} := \{ x \in \R^n\; | \; \alpha\cdot x = \beta \},   \qquad 
H_{(\alpha,\beta)}^{\le} := \{ x \in \R^n\; | \; \alpha\cdot x \le \beta \}.
$$

\noindent 
For $i\in [n]$, let $n,s_i\in \N$,
$$
M_i:=[m_{i,1},\ldots,m_{i,s_i}]^T\in \Z^{s_i\times n}, \quad 
\beta_i:=(\beta_{i,1},\ldots, \beta_{i,s_i}) \in \Z^{s_i}, 
\text{ and } \quad P_i = \{x \in \R^n \; | \; M_ix \le b_i \}.
$$ 

\noindent 
Since linear programming can be solved in polynomial-time (in the cases 
we consider) we may assume that the presentations $(n,s_i;M_i,b_i)$ are {\em 
irredundant}, i.e., $P_i$ has exactly $s_i$ facets and the sets
$P_i\cap H_{(a_{i,j},\beta_{i,j})}$, for $j\!\in\![s_i]$, are precisely the 
facets of $P_i$ for all $i\!\in\![n]$.

Now set $P:=\bigcap_{i=1}^n P_i$ and let $v\in \Q^n$. 
Note that $\size(P)$ is thus linear in $\sum^n_{i=1}\size(P_i)$. 

\begin{lem}\label{lem-NP}
{\sc Mixed-Vertex} $\in \np$.
\end{lem}

\noindent 
{\bf Proof:} 
Since the binary sizes of the coordinates of the vertices of $P$ are bounded by a polynomial in the input size, we
can use vectors $v\in \Q^n$ of polynomial size as certificates.
We can check in polynomial-time
whether such a vector $v$ is a vertex of $P$ simply by exhibiting $n$ 
facets (with linearly independent normal vectors), one from each $P_i$, 
containing $v$. If this is not the case, $v$ 
cannot be a mixed-vertex of $P$.
Otherwise, $v$ is a mixed-vertex of $P$ if and only if for each $i\!\in\![n]$ 
there exists a facet $F_i$ of $P_i$
with $v\!\in\!F_i$. Since the facets of the polytopes $P_i$ admit 
polynomial-time decriptions as $\cH$-polytopes,
this can be checked by a total of $m_1+\ldots + m_n$ polytope membership tests.

So, we can check in polynomial-time whether a given certificate $v$ is a 
mixed-vertex of $P$. Hence {\sc Mixed-Vertex} is in $\np$. \qed 

Since, in fixed dimensions we can actually list all vertices of $P$ in 
polynomial-time, one by one, it is clear that {\sc Mixed-Vertex} can be solved 
in polynomial-time when $n$ is fixed. When $n$ is allowed to vary 
we obtain hardness: 
\begin{thm}\label{thm-NP-hardness}
\label{thm:hard} 
{\sc Mixed-Vertex} is $\np$-hard.
\end{thm}

Recall that $\sqcup$ denotes disjoint union. 
The proof of Theorem \ref{thm:hard} will be based on a transformation from the 
following decision problem: 

\medskip 
\noindent 
{\sc Partition}\\ 
{\em 
Given $d\in \N$, $\alpha_1,\dots,\alpha_d\in \N$, 
is there a partition $d\!=\!I\sqcup J$ such that
$\sum_{i\in I}\alpha_i = \sum_{j\in J}\alpha_j$? \qed} 

\medskip 
Recall that {\sc Partition} was on the original list of $\np$-complete 
problems from \cite{karp72}. 

Let an instance $(d; \alpha_1,\dots,\alpha_d)$ of {\sc Partition} be given, 
and set $\alpha:=(\alpha_1,\dots,\alpha_{d})$. 
Then we are looking for a point $x\in \{-1,1\}^d$
with $\alpha\cdot x=0$.

We will now construct an equivalent instance of {\sc Mixed-Vertex}. With 
$n\!:=\!d+1$, $x:=(\xi_1,\ldots,\xi_{n-1})$ and\linebreak  
$\1_n:=(1,\ldots,1)\in \R^n$ let
$$
P_i :=\left\{\left. \begin{bmatrix} x\\ \xi_n\end{bmatrix}\right| 
-1\le \xi_i \le 1,\, -2  
\le \xi_j \le 2 \text{ for all } j\in [n]\setminus \{i\}  \right\}
$$
for $i\in [n-1]$,
$$
P_n:=\left\{\left. \begin{bmatrix} x\\ \xi_n\end{bmatrix}\; \right| 
\;  -2\cdot\1_{n-1} \le x \le 
2\cdot\1_{n-1}, \, 1\le \xi_n \le 1, \, 0 \le 2 \alpha\cdot x \le 1 \right\},
$$
and set $P\!:=\!\bigcap_{i=1}^n P_i$, 
$\widehat{\alpha}:= \begin{bmatrix}\alpha\\ 0\end{bmatrix}$. 

The next lemma shows that $P_n\cap \{-1,1\}^n$ still captures the solutions of 
the given instance of partition.

\begin{lem}\label{lem-characterization}
$(d; \alpha_1,\dots,\alpha_d)$ is a ``no''-instance of
{\sc Partition} if and only if $P_n \cap \{-1,1\}^{n}$ is empty. 
\end{lem}

\noindent 
{\bf Proof:} Suppose, first, that $(d; \alpha_1,\dots,\alpha_d)$ is a 
``no''-instance of {\sc Partition}. If $P_n$ is empty there is nothing left 
to prove. So, let $y\in P_n$ and $w\in \{-1,1\}^{n-1}\times \R$.
Since $\alpha\in \N^d$ we have $|\widehat{\alpha}\cdot w| \ge 1$. Hence, with 
the aid of the Cauchy-Schwarz inequality, we have
$$
\begin{aligned}
1 & \le |\widehat{\alpha}\cdot w| = |\widehat{\alpha}\cdot y + \widehat{\alpha}\cdot (w-y)| \le 
|\widehat{\alpha}\cdot y| + |\widehat{\alpha}\cdot (w-y)|\\
& \le \frac{1}{2}+ |\widehat{\alpha}| \cdot |w-y| = \frac{1}{2}+ |a| \cdot |w-y| 
\end{aligned}
$$ 
and thus $|w-y|\ge \frac{1}{2 |a|} > 0$. 
Therefore $P_n \cap \bigl(\{-1,1\}^{n-1}\times \R\bigr)$ is empty.

Now, let $P_n \cap \{-1,1\}^{n}=\emptyset$. Since 
$\widehat{\alpha}\in \R^{n-1}\times \{0\}$ we have
$P_n \cap \{-1,1\}^{n} =\emptyset$. \qed 

The next lemma reduces the possible mixed-vertices to the vertical edges of the standard cube.

\begin{lem}\label{lem-characterization-Pi} 
Following the preceding notation, let $v$ be a mixed-vertex of $P$. Then 
$v\!\in\!\{-1,1\}^{n-1}\times [-1,1]$.
\end{lem}

\noindent 
{\bf Proof:} 
First note that $Q\!:=\!\bigcap_{i=1}^{n-1} P_i = [-1,1]^{n-1}\times [-2,2]$. 
Therefore, for each $i\!\in\![n-1]$, the only facets of $P_i$ that meet $Q$
are those in $H_{(e_i,\pm 1)}$ and $H_{(e_n,\pm 2)}$. Since $P\subset 
[-1,1]^n$, and for each $i\in [n-1]$ the mixed-vertex $v$ 
must be contained in a facet of $P_i$, we have
$$
v \in [-1,1]^n \cap \bigcap_{i=1}^{n-1}\left( \bigcup_{\delta_i\in \{-1,1\}} H_{(e_i,\delta_i)}\right)
= \{-1,1\}^{n-1}\times [-1,1],
$$
which proves the assertion. \qed 

\medskip 
The next lemma adds $P_n$ to the consideration.

\begin{lem}\label{lem-characterization-Pn}
Let $v$ be a mixed-vertex of $P$. Then $v\in \{-1,1\}^n$.
\end{lem}

\noindent 
{\bf Proof:} 
By Lemma \ref{lem-characterization-Pi}, $v \subset \{-1,1\}^{n-1}\times 
[-1,1]$. Since the hyperplanes $H_{(e_n,\pm 2)}$ do not meet $[-1,1]^n$,
$$
v\not\in H_{(e_i,-2)} \cup H_{(e_i,2)}\quad \text{ for all } i\!\in\![n-1]. 
$$
Hence, $v$ can only be contained in the constraint hyperplanes 
$H_{(\widehat{\alpha},0)}, H_{(2\widehat{\alpha},1)}, H_{(e_n,-1)}, 
H_{(e_n,1)}$. Since $\widehat{\alpha} \in \R^{n-1}\times \{0\}$, the vector 
$\widehat{\alpha}$ is linearly dependent on $e_1,\ldots,e_{n-1}$. Hence, 
$v\!\in\!H_{(e_n,-1)} \cup H_{(e_n,1)}$, i.e., $v\in \{-1,1\}^n$. \qed 

\medskip 
Now we can prove the $\np$-hardness of {\sc Mixed-Vertex}.

\medskip 
\noindent 
{\bf Proof of Theorem \ref{thm-NP-hardness}:} 
First, let  $(d; \alpha_1,\dots,\alpha_d)$ be a ``yes''-instance of
{\sc Partition}, let $x^*:=(\xi^*_1,\ldots,\xi_{n-1}^*)\!\in\!\{-1,1\}^{n-1}$ 
be a solution, and set
$$
\xi_n^*:=1,\quad v:=\begin{bmatrix} x^*\\ \xi_n^*\end{bmatrix}, \quad
F_i:= H_{(e_i,\xi_i^*)} \cap P_i \quad \text{ for all } i\!\in\![n], 
\text{ and } \hat{F}_n:= H_{(\widehat{\alpha},0)} \cap P_n.
$$
Then $v\in \hat{F}_n \subset P_n$, hence $v\in P$ and, in fact, $v$ is a 
vertex of $P$. Furthermore, $F_i$ is a facet of $P_i$ for all $i\!\in\![n]$, 
$v\!\in\!\bigcap_{i=1}^{n} F_i$, and thus $v$ is a mixed-vertex of $P$.

Conversely, let $(d; \alpha_1,\dots,\alpha_d)$ be a ``no''-instance of
{\sc Partition}, and suppose that $v\in \R^n$ is a mixed-vertex of $P$.
By Lemma \ref{lem-characterization-Pn}, $v\in \{-1,1\}^n$.
Furthermore, $v$ lies in a facet of $P$. Hence, in particular, $v\in P_n$, 
i.e., $P_n \cap  \{-1,1\}^{n}$ is empty. 
Therefore, by Lemma \ref{lem-characterization}, 
$(d; \alpha_1,\dots,\alpha_d)$ is a ``yes''-instance of
{\sc Partition}. This contradiction shows that $P$ does not have a 
mixed-vertex.

Clearly, the transformation works in polynomial-time. \qed 

\subsection{Proof of the Second Assertion of Theorem \ref{thm:mixed}}  
We call a polyhedron {\em $\ell$-rational} if and only if it is of the form 
$\{x\!\in\!\Rn\; | \; Mx\!\leq\!b\}$ with $M\!\in\!\Q^{k\times n}$ and 
$b\!=\!(b_1,\ldots,b_k)$ satisfying $b_i\!=\!\beta_{1,i}\log|\alpha_1|
+\cdots+\beta_{k,i}\log|\alpha_k|$, with $\beta_{i,j},\alpha_j\!\in\!\Q$ for 
all $i$ and $j$. We measure the size of such a polyhedron as 
$\size(M)+\size([b_{i,j}])+\sum^k_{i=1}\size(\alpha_i)$. Clearly, it suffices 
to show that the following variant of {\sc Mixed-Vertex} is $\np$-hard: 

\smallskip
\noindent
{\sc Logarithmic-Mixed-Vertex}:\\
{\em Given $n\in \N$ and $\ell$-rational polyhedra 
$P_1,\ldots,P_n\!\subset\!\R^n$, 
does $P:=\bigcap_{i=1}^n P_i$ have a mixed vertex? \qed}

\smallskip 
Via an argument completely parallel to the last section,  
the $\np$-hardness of {\sc Logarithmic-Mixed-Vertex} follows immediately 
from the $\np$-hardness of the following variant of {\sc Partition}: 

\smallskip
\noindent
{\sc Logarithmic-Partition}\\
{\em
Given $d\in \N$, $\alpha_1,\dots,\alpha_d\in \N\setminus\{0\}$,
is there a partition $d\!=\!I\sqcup J$ such that
$\sum_{i\in I}\log\alpha_i = \sum_{j\in J}\log \alpha_j$? \qed}

\smallskip 
We measure size in {\sc Logarithmic-Partition} just as in 
the original {\sc Partition} Problem: $\sum^d_{i=1} \log \alpha_d$. 
Note that {\sc Logarithmic-Partition} is equivalent to the obvious variant 
of {\sc Partition} where we ask for a partition making the two resulting 
{\em products} be identical. The latter problem is easily seen to be 
$\np$-hard as well, via an argument mimicking the original proof of the 
$\np$-hardness of {\sc Partition} in \cite{karp72}.  \qed 

\bibliographystyle{acm}

\end{document}